\documentclass{amsart}    
\usepackage{amssymb}    
\usepackage{latexsym}    
    
\newcommand{\ZZ}{{\mathbb Z}}    
\newcommand{\RR}{{\mathbb R}}    
    
\newcommand{\NN}{{\mathbb N}}

\newtheorem{theorem}{Theorem}       
\newtheorem{remark}{Remark}       
\newtheorem{lemma}{Lemma}[section]       
\newtheorem{prop}[lemma]{Proposition}       
\newtheorem{coro}[lemma]{Corollary}       
\newtheorem{definition}[lemma]{Definition}       
\renewcommand\qedsymbol{$\Box$}    
    
\newcommand{\calw}{{\mathcal{W}}}
\newcommand{\Sub}{\mbox{Sub}}
\newcommand{\RW}{\mathcal{R}}
\newcommand{\Fbar}{\overline{F}}

\newcounter{smalllist}

\begin{document}    
\title{Uniform ergodic theorems on  subshifts over a finite alphabet}    
\author[ D.~Lenz]{Daniel Lenz}    
\maketitle    
\vspace{0.3cm}        
\noindent        
Fachbereich Mathematik, Johann Wolfgang Goethe-Universit\"at,        
60054 Frankfurt, Germany\\[2mm]  
E-mail: \mbox{ dlenz@math.uni-frankfurt.de}\\[3mm]        
\begin{abstract}    
We investigate uniform ergodic type theorems for additive and subadditive
functions on a subshift over a finite alphabet. 
We show that  every strictly ergodic subshift admits a uniform ergodic
theorem for Banach-space-valued additive functions. 
We then give a necessary and
sufficient condition on a minimal subshift to allow for a uniform subadditive
ergodic theorem.  This provides in particular a sufficient condition for
unique ergodicity. 
\end{abstract}    
    
\section{Introduction}    
Ergodic theorems for  addditive and subadditive functions play a
role in several branches of mathematics and physics. In particular, they 
are an important tool in statistical mechanics as well as in the theory of
random operators (cf. \cite{CL,GH,Rue, Sim} and references therein). 

\medskip

During recent years lattice gas models  and random operators on
aperiodic tilings  have received a lot of interest both in one
dimension and in higher dimensions (cf. \cite{BG, Dam, GH, Hof1, Hof2}
and references therein). 
In these cases one rather expects uniform ergodic theorems to hold.

\medskip

The aim of this paper is to provide a thorough study of the validity of 
such theorems in the one-dimensional case. 

In  particular, we show that
every strictly ergodic subshift over a finite alphabet admits a uniform 
additive ergodic theorem. Moreover, we give a necessary and sufficient
condition for a minimal subshift to allow for a subadditive theorem.
This gives in particular a sufficient condition for unique ergodicity. 

The  proofs are quite elementary and conceptual. Thus, it is to be
expected that a considerable  part of the material presented here,  can be extended to higher dimensional tiling dynamical systems. 

\medskip

More precisely,  we consider the following situation:\\
Let $A$ be a
finite set called the alphabet and equipped with the discrete
topology. Let $\Omega$ be a subshift over $A$. This means that $\Omega$
is a closed subset of $A^{\ZZ}$, where  $A^{\ZZ}$ is given the   
product topology and $\Omega$ is invariant under the shift operator
$T:A^{\ZZ}\longrightarrow A^{\ZZ}$, $T a(n)\equiv a(n+1)$. 

We consider
sequences over $A$ as words and  use 
standard concepts from the theory of words (\cite{Du,Lot1}).  In particular,  
$\Sub(w)$ denotes the set of subwords of $w$, the empty word is denoted
by $\epsilon$,  the number of occurrences
of $v$ in $w$ is denoted by $\sharp_v (w)$ and the length $|w|$ of the word
$w=w_1\ldots w_n$ is given by $n$.  To $\Omega$ we associate the set $\calw=\calw(\Omega)$ of
finite words associated to $\Omega$ given by $\calw\equiv
\cup_{\omega\in \Omega} \Sub(\omega)$.  A word $w\in \calw$ is called
primitive if $v$ can not be written as $v=w^l$, with $w\in\calw$ and
$l\geq 2$.  For a finite set $M$, we define $\sharp M$ to be the number
of elements in $M$.

\medskip

To phrase our additive ergodic theorem, we need the following definition.

\begin{definition}  
Let $(B,|\cdot|)$ be a Banach space. A function $F:\calw \longrightarrow  
B$ is called additive if there exists a constant $D>0$ and a  
nonincreasing function  $c:[0,\infty)\longrightarrow [0,\infty)$ with  
$\lim_{r\to \infty} c(r)=0$ s.t. the following holds  
\begin{itemize}  
\item[(A1)] $|F(v)-\sum_{j=1}^n F(v_j)|\leq \sum_{j=1}^n c(|v_j|)|v_j|$  
for $v=v_1\ldots v_n\in \calw$.   
\item[(A2)] $|F(v)|\leq D |v|$ for every $v\in \calw$.    
\end{itemize}  
\end{definition}  
Then the additive theorem  can be stated as follows.

\begin{theorem}\label{main} Let $(\Omega,T)$ be a minimal subshift over the finite
  alphabet $A$ with associated set of words $\calw=\cup_{\omega\in
    \Omega} \Sub(\omega)$. Then the following are equivalent:
\begin{itemize}
\item[{\rm (i)}] $(\Omega,T)$ is uniquely ergodic, i.e. $\lim_{|w|\to
    \infty} \frac{ \sharp_v (w)}{|w|}$ exists for all $v\in \calw.$
\item[{\rm (ii)}] The limit $\lim_{|w|\to \infty}\frac{F(w)}{|w|}$
  exists for every Banach-space-valued additive function $F$ on $\calw$.
\item[{\rm (iii)}] The limit $\lim_{n\to \infty} \frac{1}{n}\sum_{j=1}^n
  f(T^j \omega)$ exists  uniformly in $\omega\in\Omega$ for every continuous
  Banach-space-valued function  $f$ on $\Omega$.
\end{itemize}
\end{theorem}
Here, the proof of (iii) $\Longrightarrow$ (i) is well known. The proof of 
(ii) $\Longrightarrow$ (iii) follows by rather standard arguments
(cf. \cite{L2,DL}) and is in fact valid in arbitrary dimensions. So, the
hardest part is the proof of (i) $\Longrightarrow$ (ii). There,  we use
the following strategy:

\begin{itemize}
\item  A  weak form of hierarchic structure exists in arbitrary
  minimal subshifts.
\item  This structure suffices to prove a local ergodic theorem such as
  (ii). 
\end{itemize}

Here, the hierarchy found in the first step generalizes the notion of 
$n$-partition  in  Sturmian dynamical systems  investigated in 
\cite{L1,L2} (cf. \cite{DL1,DL2} as well for application to
one-dimensional quasicrystals). It is based on the concept of return
word recently introduced by Durand \cite{Du}. The second step  then uses
ideas of  \cite{GH}.

\medskip

To introduce the second result of this paper, recall that a function $F: 
\calw\longrightarrow \RR$ is called subadditive if it satisfies
$F(ab)\leq F(a)+ F(b)$.  The dynamical system $(\Omega,T)$ is said to satisfy (SET), i.e. to 
admit a uniform subadditive ergodic theorem, if, for every subadditive
function $F$, the limit
$\lim_{|w|\to \infty}\frac{F(w)}{|w|}$ exists and equals $\inf_{n\in
  \NN} F^{(n)}$, where $F^{(n)}\equiv \max\{\frac{F(v)}{|v|}: v\in
\calw, \, |v|=n\}$. 

\medskip

Define the functions $l_v:\calw\longrightarrow \RR$,  for $v\in \calw$,   and $\nu:\calw\longrightarrow \RR$  by 
$$l_v(w)\equiv (\mbox{Maximal number of disjoint copies of $v$ in
  $w$})\cdot |v|$$
and
$$ \nu(v)\equiv \liminf_{|w|\to \infty} \frac{l_v(w)}{|w|}.$$
Then a  subshift $(\Omega,T)$ over $A$ is said to satisfy
uniform positivity of quasiweights (PQ) if the following  condition holds:
\begin{itemize}
\item[(PQ)] There exists a a constant $ C>0$ with  $\nu(v) \geq C$ for all $v \in \calw$. 
\end{itemize}
Our result on subadditive ergodic theorems then reads as follows:

\begin{theorem}\label{subadditive}
Let $(\Omega,T)$ be a minimal subshift. Then the following are equivalent:  \\
(i) The subshift $(\Omega,T)$ satisfies (SET).\\
(ii) The subshift $(\Omega,T)$ satisfies  (PQ).
\end{theorem}

\begin{remark}\label{rem} {\rm  (a) As the validity of a subadditive ergodic theorem
immediately implies that the underlying subshift is uniquely ergodic, we 
see that (PQ) is a sufficient condition for unique ergodicity.\\
(b)  It is not hard to prove    (PQ) for subshifts which satisfy the
uniform positivity of weights (PW) as well as a highest power condition
(HP) given as follows:
\begin{itemize}
\item[(PW)] There exists $E>0$ with  $\liminf_{|w|\to
  \infty}\frac{\sharp_v(w)}{|w|}|v|\geq E$ for all $v\in \calw$.
\item[(HP)] There exists an $N>0$ s.t. $v^k\in\calw$ implies  $k\leq
  N$ (or equivalently: there exists a $\kappa>0$ s.t. $v$ prefix of $u
  v$ with $u$ not empty implies $|u|\geq \frac{|v|}{\kappa}$ ). 
\end{itemize}
In particular, (PQ) holds if $(\Omega,T)$ is linearly repetitive i.e. if there exists a $D>0$ s.t. every $v\in
\calw$ is a factor of every $w\in\calw$ with $|w|\geq D |v|$.  This is
further discussed in Section \ref{ex}.\\
(c) It is not hard to see that  (PQ) (or (PW)) implies in particular the
minimality of the subshift $(\Omega,T)$. In fact, minimality of
$(\Omega,T)$ is equivalent to $\nu(v)>0$ for every $v\in \calw$. From
this point of view (PQ) (or (PW)) can be seen as a very strong
minimality condition.  }
\end{remark} 
The organisation of the paper is as follows. In Section \ref{Prelim} we
review basic notions from  the theory of subshifts over finite
alphabets. The additive ergodic theorem is then contained in Section
\ref{Main}. The subadditive ergodic theorem is contained in
Section \ref{sub}. Finally, in Section \ref{ex} we provide several examples
of subshifts satisfying (PQ), thereby giving a precise form to  Remark
\ref{rem} (b).

\section{Preliminaries}\label{Prelim}
Let $\Omega$ be a subshift over the finite alphabet $A$
 with associated set $\calw$ of finite words. The following two facts
 are well known (cf. \cite{Que}):
\begin{itemize}
\item[$(*)$]  The subshift $\Omega$ is minimal, i.e. each orbit is dense if and only  
if for every $w\in \calw$ there exists an $R(w)>0$ s.t.  $w$ is a factor  
of every $v\in \calw$ with $|v|\geq R(w)$. 
\item[$(*)$] The subshift is uniquely  
ergodic if and only if the limit $\lim_{|w|\to \infty} \frac{\sharp_v  
(w)}{|w|}$ exists for every $v\in \calw$.   
  \end{itemize}

A return word of $u\in \calw$ is a word $w\in \calw$ s.t.  
$$wu\in \calw,\;\: \;\sharp_u (wu)=2,\;\:\; \mbox{$u$ is prefix of  
$wu$}.$$
 This notion was introduced by Durand in \cite{Du}.
The set of return words of $u\in \calw$ will be denoted by $\RW(u)\equiv  
\RW(u,\calw)$. Similarly, one can define for $n\in \NN$ the set  
$\RW^n\equiv \cup_{v\in \calw,\,|v|=n} \RW(v)$. If $(\Omega,T)$ is
minimal, the length of a return word of $u\in \calw$ is bounded by
$R(u)$. In particular, the set $\RW(u)$ is  
finite for every $u\in \calw$ if $\calw$ is minimal.   
  
\medskip  
  
Set $m(n)\equiv \min\{|x|: x\in \RW^n\}$. Recall that a minimal subshift  
is called aperiodic if one (and then each) of its elements is not  
periodic.
  
\begin{prop}\label{equiv}  
Let $(\Omega,T)$ be a minimal subshift. Then the following are  
equivalent:  
\begin{itemize}  
\item[{\rm (i)}] $(\Omega,T)$ is aperiodic.  
\item[{\rm (ii)}] $m(n)\longrightarrow\infty$, $n\longrightarrow  
\infty$.   
\end{itemize}  
\end{prop}  
{\it Proof.} (ii)$\Longrightarrow$(i). This is clear.\\  
(i)$\Longrightarrow$ (ii). Assume the contrary. Then, there exsists an  
$r>0$ and sequences $(u_k)$, $(v_k)$ in $\calw$ with  
$$|u_k|\leq r,\,k\in \NN; \: |v_k|\longrightarrow \infty,  
\,k\longrightarrow \infty; \:\; u_k\in \RW(v_k).$$  
As $u_k$ is a return word of $v_k$, the word $v_k$ is a prefix of $u_k  
v_k$. Thus, $v_k$ begins with $u_k^{l(k)}$, where  
$l(k)=\left[\frac{|v_k|}{|u_k|}\right]$. Here, $[x]$ denotes the largest integer not  exceeding $x$. By assumption, we have  
\begin{equation}\label{haha}  
l(k)\longrightarrow \infty, \, k\longrightarrow \infty.  
\end{equation}  
As $\Omega$ is minimal, there exists an $L>0$ s.t. every word of length  
$r$ is a factor of every word of length $L$. By \eqref{haha}, the  
word $u_k^L$ belongs to $\calw$ for sufficiently large $k$. Therefore,
the word  $u_k^L$ contains every word of length $r$ for sufficiently large  
$k$. But this implies  
$$ \sharp\{v\in \calw: |v|=r\} \leq \sharp \{v\in \Sub(u_k^L): |v|=r\}\leq |u_k|\leq r$$  
contradicting the aperiodicity of $\Omega$. \hfill \qedsymbol

\medskip

Consider now a minimal subshift $(\Omega,T)$. 
Let $u,x\in\calw$ be given. Decomposing $x$ according to the occurrences  
of $u$ in $x$ gives a unique way of writing $x$ as $x= a u_1\ldots u_l  
b$ with  
$$ a u_1\neq \epsilon, \, \sharp_u(x)=l, \, \;\mbox{$u$ prefix of  
$u_j\ldots u_l b$}, \, j=1,\ldots, l.$$  
This decomposition is called the $u$-partition of $x$. It should be  
emphasized that in minimal $\calw$ the word $u$ is a factor of every  
word of length bigger than $R(u)$ and therefore $|a|,|b|\leq R(u)$ for  
all such $w$. Note that $u_1,\ldots, u_{l-1}$ are return words of $u$  
whereas $u_l=u$. Similarly, there is a unique way of writing an arbitrary  
$\omega\in \Omega$ as $\omega=\ldots u_{-2} u_{-1} u_0 u_{1} u_{2}\dots$  
with  
\begin{itemize}  
\item[($*$)] $\omega(0)$ belongs to $u_0$,  
\item[($*$)] every occurrence of $u$ in $\omega$ begins with one of the  
$u_j$, $j\in \ZZ$.   
\end{itemize}  
If $x$ is a factor of $y$ then the $u$-partition of $y$ induces  a  
decomposition of $x$ respecting $u_1,\ldots, u_{l-1}$.   
  
\medskip  
  
For a return word $z$ of $u$ and $w\in \calw$ with $u$-partition $w=a  
u_1\ldots u_l b$, the topological number of occurrences of $z$ in $w$ is  
defined by  
$$p_{z,u}(w)\equiv \sharp\{j\in \{1,\ldots,l\}: u_j=z\}.$$  
Apparently the following proposition holds.  
  
\begin{prop} \label{frequenz} $p_{z,u}(w)=\sharp_{zu} (w)$.   
\end{prop}  
A sequence $(x_n)$ in $\calw$ is called partitioning sequence if $
|x_n|<|x_{n+1}|$.  This implies in particular  
$|x_n|\longrightarrow \infty$, $n\longrightarrow \infty$.   
  
Given a partitioning sequence $(x_n)$, one can consider the series of  
$x_n$-partitions of an $\omega\in \Omega$. This gives a hierarchic  
structure in $\omega$. This structure is a weak form  of the hierarchic  
structures present in Sturmian models or substitutional models (cf. \cite{Du,L1,L2}).  
  
\section{Uniform additive ergodic theorems}  \label{Main}
The following is the key to the proof of Theorem \ref{main}  
\begin{lemma}\label{lemmaeins}  
Let $(\Omega,T)$ be a strictly ergodic subshift over $A$ and  
let $B$ be a Banach space. Let $F:\calw \longrightarrow B$ be additive. Then  
the limit   
$\lim_{|w|\to \infty} \frac{F(w)}{|w|}$
exists.   
\end{lemma}  
{\it Proof.} We consider two cases:\\
{\it Case 1.} $(\Omega,T)$ is aperiodic.\\
Let $(x_n)$ be an arbitrary partitioning sequence. As
$(\Omega,T)$ is minimal, there exists only finitely many return words
for each  $x_n$, $n\in \NN$. Therefore,    
the return words of $x_n$ can be denoted by  $y^{(n)}_1,\ldots
y^{(n)}_{l(n)}$ with a suitable $l(n)\in \NN$. As  
$(\Omega,T)$ is strictly ergodic, the limits  
\begin{equation}\label{oho} \nu_j^{(n)}\equiv \lim_{|w|\to \infty} \frac{ p_{y^{(n)}_j, x_n}(w)}{|w|}
  |y^{(n)}_j|= \lim_{|w|\to \infty} \frac{ \sharp_{y^{(n)}_j
    x_n}(w) }{|w|} |y^{(n)}_j|
\end{equation}  
exists by Proposition \ref{frequenz} for all $j,n\in \NN$ with $j\leq  
l(n)$. Set  
$$ F^{(n)}\equiv \sum_{j=1}^{l(n)} \nu_j^{(n)}  \frac{F(y^{(n)}_j)}{|y^{(n)}_j|}.$$  
Of course, it suffices to show that for every $\varepsilon>0$ there  
exists an $n\in \NN$ with   
$ |\frac{F(w)}{|w|}- F^{(n)}|\leq \varepsilon$ for all $w\in \calw$ long  
enough.  So, choose $\varepsilon>0$ arbitrary. By Proposition  
\ref{equiv} and the additivity of $F$, there exists an $n\in \NN$ with  
\begin{equation}\label{huhu}  
c( m(|x_n|))\leq \frac{\varepsilon}{2}.  
\end{equation}  
Considering now a $w\in \calw$ with $x_n$-partition $w=a u_1 \ldots u_s  
b$, we can estimate  
{\footnotesize
\begin{eqnarray*}  
& &\left|\frac{F(w)}{|w|} - F^{(n)}\right|=  \left|\frac{F(w)}{|w|}- \sum_{j=1}^{l(n)}  
\nu_j^{(n)}  \frac{ F(y^{(n)}_j)}{|y^{(n)}_j|} \right|\\  
&\leq&\left| \frac{F(w)}{|w|} - \frac{ F(a)+ \sum_{j=1}^s   F(u_j) +
    F(b)}{|w|}\right| +\left| \frac{ F(a)+ \sum_{j=1}^s   F(u_j) +F(b)}{|w|} - \sum_{j=1}^{l(n)}  \nu_j^{(n)}  \frac{ F(y^{(n)}_j)}{|y^{(n)}_j|}\right|\\
&\leq& \frac{ c(|a|)|a| + c(|b|)|b|}{|w|} + \sum_{j=1}^s \frac{c(|u_j|)}{|w|}
  |u_j| +  \frac{ D |a| + D|b|}{|w|}\\
 &+& \sum_{j=1}^{l(n)}  \left| \frac{
      p_{y^{(n)}_j, x_n}(w)  |y^{(n)}_j| }{|w|}- \nu_j^{(n)}\right|
  \frac{|F(y^{(n)}_j)|}{|y^{(n)}_j|}\\
&\leq & \frac{ (c(0) + c(0)) R(x_n)}{|w|}+ \frac{2 D R(x_n)}{|w|} +
  \frac{\varepsilon}{2} + D  \sum_{j=1}^{l(n)} \left| \frac{p_{y^{(n)}_j, x_n}(w)  |y^{(n)}_j|}{|w|}- \nu_j^{(n)}\right|,
\end{eqnarray*}  }
where we used \eqref{huhu} in the last equation. By \eqref{oho}, this gives
$ \left|\frac{F(w)}{|w|}- F^{(n)}\right|\leq \varepsilon$
for all $w$ which are sufficiently long and the  proof of Case 1 is
finished.\\
\\{\it Case 2.} $(\Omega,T)$ is periodic. \\
There exists a $p\in \NN$ with $T^p \omega=\omega$ for every $\omega\in
\Omega$. Fix a word $v$  in $\calw$ of length $p$. Let
$n\in \NN$ be arbitrary. Then every $w\in \calw$ with $|w|\geq 3
|v^n|$ can be written as $w= a v_1 \ldots v_s b$ with $v_j=v^n$,
$j=1,\ldots, s,$ and $|a|,|b|\leq |v^n|$. Using this partition of $w$
instead of the $x_n$-partition, one can mimick the proof of Case 1. This 
gives the existence of the limit in question.   \hfill \qedsymbol

\begin{lemma}\label{lemmazwei} Let $\Omega, \calw$ and $B$ be as in the
  foregoing lemma. Let
  $f:\Omega\longrightarrow B$ be a continuous function. Then,
  $\frac{1}{n} \sum_{j=1}^n f(T^j \omega)$ converges uniformly in
  $\omega$ towards a constant function. 
\end{lemma}
{\it Proof.} To $\omega\in \Omega$, we can associate the function
$F_\omega:\calw \longrightarrow B$ defined by 
$$ F_\omega(w)\equiv \sum_{j=0}^{|w|-1} f( T^{k_\omega^w+j} \omega),$$
where $k_\omega^w\in \NN$ is minimal with $\omega_{k_\omega^w } \ldots
\omega_{k_\omega^w+   |w|-1}=w$. 

\medskip

{\bf Claim.} There exists a nonincreasing function
$c:[0,\infty)\longrightarrow [0,\infty)$ with $\lim_{r\to \infty}
c(r)=0$ s.t. 
$$ |F_\omega(w)- F_\rho(w)|\leq c(|w|) |w|.$$
{\it Proof of the claim.} It is clearly enough to show that, for every
$\varepsilon >0$, there exists an $l=l(\varepsilon)$ with $|F_\omega(w)-
F_\rho(w)|\leq \varepsilon |w|$ for all $\omega,\rho\in \Omega$ and
every $w\in \calw$ with $|w|\geq l$. \\
As $f$ is continuous on the compact space $\Omega$, there exists a $k\in
\NN$ s.t. $\omega_{-k}\ldots \omega_{k}= \rho_{-k}\ldots \rho_k$ implies 
\begin{equation} \label{hoho}
|f(\omega)-f(\rho)|\leq \frac{\varepsilon}{2}.
\end{equation}
Set $l\equiv \max\{\frac{ 8 k \|f\|_\infty}{\varepsilon}, 2k\}$, where
$\|\cdot\|_{\infty}$ denotes the supremum-norm. For $w\in \calw$
with $|w|\geq l$, we can estimate the  value of $D\equiv
|F_\omega(w)-F_\rho(w)|$ by 
{\small
\begin{eqnarray*}
 |\sum_{j=0}^{|w|-1} (f(T^{k_\omega^w +j}
\omega)- f(T^{k_\rho^w +j} \rho)|
&\leq & 4 \|f\|_\infty k + \sum_{j=k}^{|w|-k-1} | (f(T^{k_\omega^w +j}
\omega)- f(T^{k_\rho^w +j} \rho)|\\
\eqref{hoho}\,\leq 4 \|f\|_\infty k + (|w|-2 k)\frac{\varepsilon}{2}
&=& |w|\left( \frac{ 4 \|f\|_\infty k}{|w|} + \frac{(|w|-2 k)}{|w|}
  \frac{\varepsilon}{2}\right) \leq |w| \varepsilon.
\end{eqnarray*}}
The proof of the claim is finished. 

\medskip

The claim implies that
\begin{itemize}
\item[$(*)$] for every $\omega\in \Omega$ the function $F_\omega$ is additive,
\item[$(*)$]  $\left\| \frac{ F_\omega(w)-F_\rho(w)}{|w|}\right\| \leq
  c(|w|)$, with $\lim_{|w|\to \infty} c(|w|)=0$. 
\end{itemize}
Using this and Lemma \ref{lemmaeins}, it is straightforward to finish the proof of the lemma. \hfill \qedsymbol

\medskip

Now, the proof  of Theorem \ref{main} can easily be accomplished. 

\medskip

{\it Proof of Theorem \ref{main}}. The implication (i) $\Longrightarrow$
(ii) follows immediately from Lemma \ref{lemmaeins}. Similarly,
(ii) $\Longrightarrow$ (iii) follows immediately from Lemma
\ref{lemmazwei}. Finally, (iii)$\Longrightarrow $ (i) is clear as the
function $w\mapsto \sharp_v (w)$ is additive for each $v\in \calw$. The
proof of the theorem is finished. \hfill \qedsymbol

\section{Uniform subadditive ergodic theorems}\label{sub}
Here, we consider subadditive functions on $\calw$. The key results are the 
following two  lemmas. Their proofs use and considerably extend  ideas
from \cite{L1,L2}.

\begin{lemma} Let $(\Omega,T)$ be a minimal   subshift over $A$ satisfying
  (PQ). Then, the subshift  $(\Omega,T)$ satisfies  (SET) as well. 
\end{lemma}
{\it Proof.} Set $\Fbar\equiv  \inf_{n\in \NN} F^{(n)}$. We will show
$$ (A)\;\:  \limsup_{|w|\to \infty} \frac{F(w)}{|w|}\leq
\Fbar\;\:\mbox{ and}\; (B)\;\: \liminf_{|w|\to \infty} \frac{F(w)}{|w|}\geq \Fbar.$$

Ad (A): It is clearly enough to show $\limsup_{|w|\to \infty}
\frac{F(w)}{|w|}\leq F^{(n)}$ for every $n\in  \NN$. But this follows
  easily form the subadditivity of $F$ by chopping each  $w$  into parts
  of length $n$ and a boundary word  and using that the boundary terms
  tend to zero (cf. \cite{L1,L2}). \\
Ad (B): Assume the contrary. This implies, in particular, $\Fbar>-\infty$ 
and that there exists a sequence $(v_n)$ in
$\calw$ as well as a  $\delta>0$ with $|v_n|$ tending to $\infty$ for $n\to
\infty$ and 
\begin{equation}\label{oo}
 \frac{F(v_n)}{|v_n|} \leq \Fbar -\delta
\end{equation}
for  every $n\in \NN$. 
Moreover, by (A), there exists  an $L_0\in \RR$ with
\begin{equation}\label{uu}
 \frac{F(w)}{|w|} \leq \Fbar + \frac{ C \delta}{8}
\end{equation}
for all $w\in \calw$, $|w|\geq L_0$, where $C$  is the constant from
(PQ). 

\medskip

Fix  $m\in\NN$ with $|v_m|\geq L_0$. Using (PQ), we can now find an
$L_1\in\RR$ s.t. every $w\in \calw$ with $|w|\geq L_1$ can be written  as $w=  x_1  v_m x_2 v_m \ldots x_l v_m x_{l+1}$ with
\begin{equation}\label{hatata}
\frac{l-2}{2} \geq \frac{C}{4} \frac{|w|}{|v_m|}.
\end{equation} 
Now, considering only every other copy of $v_m$ in $w$, we can write $w$ as
$ w=y_1 v_m y_2 \ldots y_r v_m y_{r+1}$, with $|y_j|\geq |v_m|\geq L_0$, 
$j=1,\ldots,r+1,$ and by \eqref{hatata} 
$$r\geq \frac{l-2}{2} \geq \frac{C}{4} \frac{ |w|}{ |v_m|}.$$
Using  $\eqref{oo}$, $\eqref{uu}$ and  this estimate, we can now calculate
\begin{eqnarray*}
\frac{F(w)}{|w|} &\leq & \sum_{j=1}^{r+1} \frac{ F(y_j)}{|y_j|}
\frac{|y_j|}{|w|} +  \frac{F(v_m)}{|v_m|}\frac{ r |v_m|}{|w|}
\leq  \sum_{j=1}^{r+1} (\Fbar +\frac{C}{8}\delta )\frac{|y_j|}{|w|} + (
\Fbar -\delta) \frac{ r |v_m|}{|w|}\\
&\leq& \Fbar + \frac{C}{8}\delta    - \frac{C}{4} \frac{ |w|}{ |v_m|}\frac{|v_m|}{|w|}
\delta 
\leq \Fbar - \frac{C}{8}\delta. 
\end{eqnarray*}

As this holds for arbitrary $w\in \calw$ with $|w|=L_1$, we arrive at
the obvious contradiction $F^{(L_1)}\leq \Fbar - \frac{C}{8}\delta
< \inf_{n\in \NN} F^{(n)}$. This finishes the proof.  \hfill \qedsymbol

\medskip

\begin{lemma} Let $(\Omega,T)$ be a minimal subshift over $A$ satisfying 
  (SET). Then,  the subshift  $(\Omega,T)$ satisfies (PQ) as well. 
\end{lemma}
{\it Proof.} Note that, for $v\in \calw$,  the function $(-l_v)$ is
subadditive. Thus, the equation 
\begin{equation}\label{lalala}
\nu(v)\equiv \liminf_{|w|\to \infty} \frac{l_v(w)}{|w|}=\lim_{|w|\to \infty}
\frac{l_v(w)}{|w|}
\end{equation}
holds by (SET). The proof will now  be given by
contraposition. So, let us  assume that the values $\nu(v)$, $v\in
\calw$,  are not bounded away from zero.  
As the system is minimal, we have
$\nu(w)>0$ for every $w\in \calw$. Thus, there exists a sequence
$(v_n)$ in $\calw$ with 
\begin{equation}\label{uuu}\nu(v_n)>0, \;\:\mbox{and}\;\:    \nu(v_n)\longrightarrow 0, \;\:\; n\to \infty.
\end{equation}
As the alphabet $A$ is finite, there are only finitely many words of a
prescribed length. Thus, \eqref{uuu} implies
\begin{equation}\label{vvv} |v_n|\longrightarrow \infty, \;\:\; n\to \infty.
\end{equation}
Replacing $(v_n)$  by a suitable subsequence,  we can assume by
\eqref{uuu} that the equation  
\begin{equation}\label{cucuc}
\sum_{n=1}^\infty \nu(v_n)<\frac{1}{2}
\end{equation}
holds.  
By $\eqref{lalala}$,  $(\ref{vvv})$ and $\eqref{cucuc}$, we can choose inductively for each  $k\in
\NN$ a number $n(k)$, with
\begin{equation}\label{ccc}
\sum_{j=1}^k \frac{l_{n(j)}(w)}{|w|} <\frac{1}{2}
\end{equation}
 for every 
  $w\in \calw$ with $|w|\geq\frac{|v_{n(k+1)}|}{2}$. Note that
  $\eqref{ccc}$ implies
\begin{equation}\label{blablabla}
|v_{n(k)}|< \frac{|v_{n(k)+1}|}{2}
\end{equation}
as $\frac{l_{n(k)}(v_{n(k)})}{|v_{n(k)}|}=1$. 
Define the function  $l:\calw\longrightarrow \RR$ by 
$$l(w)\equiv \sum_{j=1}^{\infty} l_{n(j)}(w).$$
Note that the sum is actually finite for each $w\in \calw$. Obviously, $(-l)$ is subadditive. Thus, by assumption, the limit
$\lim_{|w|\to \infty} \frac{l(w)}{|w|}$
exists. On the other hand, we clearly have
$$\frac{l(v_{n(k)})}{|v_{n(k)}|}\geq
\frac{l_{n(k)}(v_{n(k)})}{|v_{n(k)}|}\geq 1$$
as well as by the induction construction
$\eqref{ccc},\eqref{blablabla}$ 
$$\frac{l(w)}{|w|}=\sum_{j=1}^k \frac{l_{n(j)}(w) }{|w|}<\frac{1}{2}$$
for $w\in \calw$ with $\frac{ |v_{n(k+1)}| }{2} \leq |w|< |v_{n(k+1)}|$.
This gives a contradiction proving the lemma.\hfill $\Box$\\

\medskip

{\it Proof of Theorem \ref{subadditive}.} This theorem follows
immediately from the foregoing two lemmas. \hfill \qedsymbol

\section{Examples}\label{ex}
In this section we discuss two classes of  examples satisfying condition
(PQ). We will need the following proposition proved similarly to Proposition
\ref{equiv}. 

\begin{prop} \label{disjoint} For $(\Omega,T)$ the following are
  equivalent:\\
{\rm (i)} $\Omega$ satisfies (HP).\\
{\rm (ii)} There exists a $\kappa>0$, s.t. the length of every return word of
$v$ is at least $\frac{|v|}{\kappa}$.
\end{prop}
{\it Proof.} (i) $\Longrightarrow$ (ii).  Let $u$ be a return word to $v$. Then $v$ is a prefix of
$uv$. Thus, $v$ starts with $u^l$, where $l=[ \frac{|v|}{|u|} ]$. By (HP)
this implies $[\frac{|v|}{|u|}]\leq N$ which in turn  yields
$\frac{|v|}{|u|}-1\leq N$. Now, (ii) follows easily. \\
(ii) $\Longrightarrow$ (i). Let $v\in \calw $ be primitive  and assume that $v^l$ belongs to $\calw$ as well for a suitable
$l\in \NN$, $l\geq 2$. Then $v$ is a
return word to $v^{l-1}$ implying $|v|\geq\frac{|v^{l-1}|}{\kappa}= \frac{l-1}{\kappa} |v|$.  This gives
immediately $\kappa \geq l-1$ and the proof of the proposition is finished.   
 \hfill \qedsymbol

\medskip

From this proposition, we can derive a   sufficient
condition for (PQ). 

\begin{prop} If $(\Omega,T)$ satisfies (HP) and (PW), then it satisfies
  (PQ) as well and (SET) holds. 
\end{prop}
{\it Proof. } Choose an arbitrary $v\in \calw$. By (PW), a word $w\in
\calw$ of sufficient length contains at least $\frac{E}{2}
\frac{|w|}{|v|}$ copies of $v$, where $E$ is the constant from (PW). By
(HP) and Proposition \ref{disjoint} (ii), this implies that a word $w$
of sufficient length contains at least $\frac{E}{2}
\frac{|w|}{|v| (\kappa +1) }$ disjoint copies of $v$. This gives
$\nu(v)\geq \frac{E}{2(\kappa +1)}>0$. Thus, $(\Omega,T)$ satisfies
(PQ). Moreover, (PQ) implies minimality of $(\Omega,T)$  (cf. Remark 1. (c)). Now, the
proposition follows from Theorem \ref{subadditive}.  \hfill \qedsymbol. 

\begin{coro} If $(\Omega,T)$ is linearly repetitive, then a subaddditive
  ergodic theorem holds. 
\end{coro}
{\it Proof.} It is well known (and easy to see) that a linearly
repetitive $(\Omega,T)$ satisfies (HP) and (PW). Thus, the corollary
follows  from the foregoing proposition. \hfill \qedsymbol

\medskip
\begin{remark}{\rm  (a)  In fact, as shown in \cite{L2},  linear repetitivity implies a subadditive
ergodic theorem in  tiling dynamical systems of arbitrary dimension (cf.
\cite{DL} as well). \\
(b) The relationship between linear repetitivity and (HP) +(PW) is
currently under investigation in \cite{DL4}}.
\end{remark}

\medskip

{\it Acknowledgements.} The author would like to thank David Damanik for 
many stimulating discussions. In fact, this paper would not have been
possible without the collaboration leading to  \cite{DL,DL4}.

\end{document}